\documentclass[production,11pt]{amsart}
\usepackage{amssymb}
\usepackage{amsmath,amscd}
\usepackage{combelow}
\usepackage{mathrsfs}

\newtheorem*{lemma*}{Lemma}
\newtheorem*{theorem*}{Theorem}
\newtheorem*{remark*}{Remark}

\theoremstyle{remark}

\theoremstyle{definition}

\newtheorem*{definition*}{Definition}

\newcommand{\R}         {{\mathbb{R}}}
\newcommand{\rmap}{\longrightarrow}
\newcommand{\diffto}{\xrightarrow{\raisebox{-0.2 em}[0pt][0pt]{\smash{\ensuremath{\sim}}}}}

\providecommand{\rest}[2]{\left. #1 \right|_{#2}}

\begin{document}

\title[On cohom. obstructions for log symplectic struct.]{On cohomological obstructions for the existence of log symplectic structures}
\author[I. M\u{a}rcu\cb{t}]{Ioan M\u{a}rcu\cb{t}}
\address{Depart. of Math., Utrecht University, 3508 TA Utrecht, The Netherlands}
\email{i.t.marcut@uu.nl}
\author[B. Osorno Torres]{Boris Osorno Torres}
\address{Depart. of Math., Utrecht University, 3508 TA, Utrecht, The Netherlands}
\email{b.osornotorres@uu.nl}

\begin{abstract}
We prove that a compact log symplectic manifold has a class in the second cohomology group whose powers, except
maybe for the top, are nontrivial. This result gives cohomological obstructions for the existence of log symplectic
structures similar to those in symplectic geometry.
\end{abstract}

\maketitle

A Poisson structure $\pi$ on a smooth manifold $M$ of dimension $2n$ is called a \emph{log symplectic} structure
if the map
\[\wedge^n \pi: M \rmap \bigwedge\nolimits^{\!2n} TM, \ \ x \mapsto \wedge^n \pi(x)\]
is transverse to the zero section.

These structures were initially studied in the framework of deformation quantization in \cite{Nest.1996}, where they are called \emph{b-symplectic} structures. Later, their complex analogue was considered in \cite{Goto.2000}, where they were first given the name log symplectic. In the context of Poisson geometry, this class of Poisson structures was introduced on two-dimensional surfaces in \cite{Radko.2002} (under the name of \emph{topologically stable Poisson structures}) where a complete classification was obtained. In higher dimensions a systematic investigation of the geometric properties of log symplectic structures appeared
in \cite{Guillemin.Eva.2012}. Their integrations by symplectic groupoids were studied in \cite{Gualtieri.2012}.

Our interest in log symplectic structures comes from the fact that these can be used to construct regular corank-one
Poisson structures. First, the singular locus of a log symplectic structure $Z:=(\wedge^n\pi)^{-1}(0)$ (if nonempty)
carries a regular corank-one Poisson structure with a very special property: it has a transverse Poisson vector field
\cite{Guillemin.Eva.2012}. Secondly, a log structure can be used to construct a regular corank-one Poisson on
$M\times S^1$, simply given by
\[\pi+X\wedge \frac{\partial}{\partial \theta},\]
where $X$ is the modular vector field of $(M,\pi)$. However, our result excludes the possibility of using this procedure
to construct corank-one Poisson structures in some interesting examples, e.g.\ on $S^4\times S^1$.

Our result is the following:
\begin{theorem*}
Let $(M^{2n},\pi)$ be a compact log symplectic manifold. Then there exists a class $c\in H^2(M)$ such that
$c^{n-1}\in H^{2n-2}(M)$ is nonzero.
\end{theorem*}

\begin{proof}
Denote by $Z:=(\wedge^n\pi)^{-1}(0)$ the singular locus of $\pi$. If $Z=\emptyset$, we can apply the usual argument
from symplectic geometry. Assume that $Z\neq \emptyset$.

We first assume that $M$ is orientable. Let $\mu$ be a volume form on $M$ and denote by
$t:=\langle\pi^n,\mu\rangle$. The singular locus becomes $Z=\{t=0\}$. The log condition implies that $t$ is a
submersion along $Z$, so we can find a retraction $r:U\to Z$, where $U$ is an open around $Z$, such that
$(r,t):U\diffto Z\times (-\delta,\delta)$ is a diffeomorphism. Since $Z$ is a Poisson submanifold (it is fixed by all Poisson
automorphisms, hence all Hamiltonians are tangent to $Z$), in this open we can write $\pi=t\partial/\partial t\wedge
X_t+w_t$ for a vector field $X_t$ and a bivector $w_t$ on $Z$, both depending smoothly on $t$. Since
$1/t\pi^n=n\partial/{\partial t}\wedge X_t\wedge w_t^{n-1}$ is nowhere vanishing, we have that the bivector
$\partial/\partial t\wedge X_t+w_t$ is invertible. Denote its inverse by $\alpha_t\wedge dt +\beta_t$, with $\alpha_t$
and $\beta_t$ forms on $Z$ depending smoothly on $t\in (-\delta,\delta)$. Then $\omega:=\pi|_{M\backslash Z}^{-1}$ can be written as
\[\rest{\omega}{U\backslash Z}=\alpha_t\wedge dt/t+\beta_t.\]
Since $\omega$ is closed we get that $\alpha_0$ and $\beta_0$ are closed, and since $dt\wedge \alpha_0+\beta_0$ is
invertible, it follows that $\alpha_0\wedge\beta_0^{n-1}$ is a volume form on $Z$. Since $Z$ is compact, this implies
that $\beta_0^{n-1}$ cannot be exact. We will construct a closed 2-form $\omega'$ on $M$ whose pullback to $Z$ is
$\beta_0$; hence $c:=[\omega']$ will satisfy the conclusion of the theorem.

Let $\chi:(-\delta,\delta)\rightarrow \R$ be a bump function that takes the value 1 for $|t|\leq \delta/4$, and 0 for $|t|\geq
\delta/2$. Consider the 2-form $\omega'$ on $M\backslash Z$ that coincides with $\omega$ outside of $U$ and on
$U\backslash Z$ it is given by
\[\omega'|_{U\backslash Z}=(\alpha_t-\chi(t)\alpha_0)\wedge dt/t+\beta_t.\]
$\omega'$ extends smoothly to $Z$, since for $|t|\leq \delta/4$ it can be written as $\omega'=\lambda_t\wedge dt
+\beta_t$, where $\lambda_t=\int_0^1\dot{\alpha}_{ts}ds$, or equivalently $\alpha_t=\alpha_0+t\lambda_t$. So
$\omega'$ is a closed 2-form on $M$ whose pullback to $Z$ is $\beta_0$; thus $[\omega']^{n-1}\neq 0$.

If $M$ is not orientable, consider $p:\widetilde{M}\to M$ the orientable double cover, and let
$\gamma:\widetilde{M}\diffto \widetilde{M}$ be the corresponding deck transformation. We first construct a tubular
neighborhood $(\widetilde{r},t):\widetilde{U}\diffto \widetilde{Z}\times (-\delta,\delta)$ of the singular locus
$\widetilde{Z}:=p^{-1}(Z)$ of $\widetilde{\pi}:=p^*(\pi)$, with $\widetilde{U}=p^{-1}(U)$, and such that the action
of $\gamma$ corresponds to $\gamma(z,t)=(\gamma(z),-t)$, for $(z,t)\in \widetilde{Z}\times (-\delta,\delta)$. The map
$\widetilde{r}:\widetilde{U}\to \widetilde{Z}$ can be constructed by lifting a retraction $r:U\to Z$. Consider a volume
form $\mu_0$, and denote by $f$ the smooth function satisfying $\gamma^*(\mu_0)=-e^{f}\mu_0$. Then the volume
form $\mu:=e^{f/2}\mu_0$ satisfies $\gamma^*(\mu)=-\mu$. Thus, by shrinking $U$, we can use
$t:=\langle\widetilde{\pi}^n,\mu\rangle$ to construct the desired tubular neighborhood. As before, on
$\widetilde{Z}\times (-\delta,\delta)$ we can write $p^*(\omega|_{U\backslash Z})=\alpha_t\wedge dt/t +\beta_t$.
Invariance under $\gamma$ implies that $(\gamma|_{\widetilde{Z}})^*(\alpha_t)=\alpha_{-t}$ and
$(\gamma|_{\widetilde{Z}})^*(\beta_t)=\beta_{-t}$. In particular $\alpha_0$ and $\beta_0$ are invariant. Thus,
choosing the function $\chi(t)$ from the construction from the orientable case to satisfy $\chi(t)=\chi(-t)$, we obtain an
invariant closed 2-form $\omega'$ on $\widetilde{M}$ that satisfies $[\omega']^{n-1}\neq 0$. Invariance implies that
$\omega'=p^*(\omega'')$ for a closed 2-form $\omega''$ on $M$; hence $c:=[\omega'']$ satisfies the conclusion.
\end{proof}

\begin{remark*}\rm
Observe that for $Z\neq\emptyset$ the proof of the theorem uses only the compactness of $Z$ and not that of $M$.
\end{remark*}

\noindent \textbf{Acknowledgements.} We would like to thank M. Crainic for useful discussions. The first author was
supported by the ERC Starting Grant no. 279729 and the second by the NWO VIDI project ``Poisson Topology'' no.
639.032.712.

\end{document}